\documentclass[12pt]{article}
\usepackage{txfonts}
\usepackage{graphicx}
\title{\Large \bf \boldmath\ \\ Suita Conjecture for a Complex Torus} 
\author{\large  Robert Xin DONG$^\ast$} 
\date{}
\begin{document}
\maketitle

\renewcommand{\thefootnote}{\fnsymbol{footnote}}

% µ¥Î»¡¢µØÖ·¡¢»ù½?
\footnotetext{\hspace*{-5mm} \begin{tabular}{@{}r@{}p{15.4cm}@{}}
$^\ast$ & Graduate School of Mathematics, Nagoya University, Nagoya 464-8602, Japan \& Department of Mathematics, Tongji University, Shanghai 200092, China. \\
& Email: 1987xindong@tongji.edu.cn\\
& Received September 10, 2012; in final form, January 10, 2013\\
\end{tabular}}

\renewcommand{\thefootnote}{\arabic{footnote}}

\begin{abstract} % ÕªÒª
\fontsize{11}{12}
\selectfont
  \,The author proves that the generalized Suita conjecture holds for any complex torus, which means that $ \alpha\pi K \geq c^2(\alpha\in\mathbb R)$,  $c$ being the modified logarithmic capacity and $K$ being the Bergman kernel on the diagonal. The open problems for general compact Riemann surfaces with genus $\geq2$ is also elaborated. The proof relies in part on elliptic function theories.\\
  

\end{abstract}

\section{Introduction}
In complex analysis in one variable, a famous question is to study the relations between two functions: the Bergman kernel and the logarithmic capacity. Such a problem has been well studied since 1970s by a number of mathematicians and it even requires tools in several complex variables. Our interest is mainly on compact Riemann surface case, and in this paper we first deal with a complex torus.

Our main theorem is stated as follows:

\textbf{Theorem 1.1\,\,\,\,(i)}  For any complex torus

 $$X_\tau:=\mathbb C/(\mathbb Z+\tau\mathbb Z)\, \,\, \,\, \,\,\big(\tau\in \mathbb C \,\,  and   \,\, Im \,\tau>0\big),$$ we know that

$$\alpha\,\pi\, K \geq c^2,\,\,\,\,\alpha \approx 6.2034, \,\,\ $$
$c$ being the modified logarithmic capacity, and $K$ being the Bergman kernel on the diagonal; \\

%\emph{(here $c$ being the modified logarithmic capacity and $K$ being the Bergman kernel on the diagonal);}

\textbf{\hspace{2.4cm} (ii)} {  $``="$ is attainable when $Im \tau \approx 1.9192.$\\}

This theorem can be seen as a generalized version of Suita conjecture, which we will explain later in detail. And as a special case when $$Im \tau \to +\infty,$$ the following corollary holds.

\textbf{Corollary 1.1.}\, \,\,\, { In particular, we know
$$\lim_{Im \tau \to +\infty}\,\frac{\pi K}{c^2}=+\infty.$$}

%\hspace{0.4cm}

For the open problem for compact Riemann surfaces with genus $\geq2$, we will elaborare in the last part of this paper and it would be a fruitful task to achieve new results about that case.

\section{Preliminaries} \label{section2}

\subsection {Suita conjecture for a bounded in $\mathbb C.$}

Let us first recall some basic notations and facts. The (negative) Green function for a bounded domain $D\subset \mathbb C$ satisfies the following:
\begin{displaymath}
\left\{ \begin{array}{ll}
\Delta G_D(\cdot, z)=2\pi \delta_z\\
G_D(\cdot, z)=0\ \ \,\,on  \ \partial D.\\
\end{array} \right.
\end{displaymath}

Define $c_D(z):=\exp  \lim_{\zeta \to z}(G_D(\zeta, z)-\log|\zeta-z|)$ to be the logarithmic  capacity of $\mathbb C\setminus  D$ with respect to $z$. On the other hand, we know that the Bergman kernel on the diagnal is $K_D(z):=\sup\,\{|f(z)|^2: f \,\, holomorphic\,\,  in\,\, D, \,\,{\int_D {|f|^2 d\lambda}}\leq 1\}$. \\

In 1972, Suita [1] conjectured  that: \, $c^2_D \leq\,  \pi \,K_D,$ \,\,\,$\forall z \in D.$\,\,\,\,(Suita proposed this conjecture on open Riemann surfaces admitting Green functions.)

Such a conjecture has an equivalent geometric interpretation as follows: $$Curv\, {c_D|dz|}\leq -4,$$ since

$$K_D=\frac{1}{\pi}\frac{\partial^2}{\partial z\partial \bar z}(\log c_D)\,.$$

For some cases, such as a simply connected domain, the $``="$ can be achieved; and if the domain is an annulus, then $`` < "$ always holds. However, for other cases, it seems so difficult to give an answer to this conjecture by direct computations, that we may need a tool.

\subsection {Ohsawa-Takegoshi $L^2$ extension theorem}

The relations between Suita conjecture and $L^2$ Extension Theorem was first observed by Ohsawa[2], who proved that $\,c^2_D \leq\, 750\, \pi \,K_D,$ \,$\forall z \in D.$ \,Later this constant on the right-hand side was improved by a number of authors (see [3-7]).
%Siu$^{[3]}$,Berndtsson$^{[4]}$,Chen$^{[5]}$, B\l{}ocki$^{[6]}$ and Guan-Zhou-Zhu$^{[7]}$.\\

In early 2012, B\l{}ocki obtained the following beautiful result (see[8]).

$\textbf{Theorem 2.1\ \ } $(Ohsawa-Takegoshi $L^2$ extension theorem with optimal constant) \\
{ Let $\Omega$ be a pseudoconvex  domain in $\mathbb C ^{n-1}\times D,$\, where $D$ is a bounded domain in $\mathbb C$ containing the origin. Then for any holomorphic $f$ in $\Omega^ {\prime}:=\Omega \bigcap \{z_n=0\}$ and $\varphi$ plurisubharmonic in $\Omega,$ one can find a holomorphic extension $F$ of $f$ to $\Omega$  such that}
$$\int_\Omega |F|^2e^{-\varphi}d\lambda
\leq \frac {\pi}{(c_D(0))^2}\int_{\Omega^{\prime}} |f|^2e^{-\varphi}d\lambda^{\prime}.$$

This means Suita conjecture holds for any bounded domain $D$ in the complex plane $\mathbb C$. And during recent months Guan-Zhou[9] have made great progress by showing this conjecture, originally stated for open Riemann surfaces admitting Green functions, is true. Thus, it would be very interesting to generalize similar results to the compact Riemann surfaces case.

\subsection {The case of a complex torus}

In this paper, we treat only the compact case, whose Green function(in the usual sense) does not even exist! But after modifying the definition, we can still consider similar problems. Meanwhile, since compact Riemann surfaces can be classified by the genus, and the case for a Riemann sphere is trivial. 

Therefore, we will first deal with a complex torus, which is denoted by $X_\tau:=\mathbb C/\big(\mathbb Z+\tau\mathbb Z\big)\,\big(\tau\in \mathbb C$, Im$\,\tau>0\big).$

 Let $H^{1,0}_{(2)}(X)$ be the Hilbert space of holomorphic 1-forms on $X_\tau$ such that $\big|\int_{X_\tau} f\wedge \bar f\,\big| <\infty.$ According to the definition, we know that the Bergman kernel of $X_\tau$ on the diagonal is the (1,1)-form $$K_{X_\tau}(z)=\frac{1}{Im\,\tau} dz\wedge d\bar z.$$

Then we consider a function $g(z,w):  {X_\tau}\times {X_\tau} \to [-\infty,0),$  s.t. for each fixed $w \in {X_\tau}:$ \\
(a)\ \ $\Delta_{\omega} g(\cdot,w)=-1$   on ${X_\tau}\setminus \{w\}$\, where $\Delta_{\omega}$ is the Laplacian with respect to the metric $\omega= \frac{1}{Im\,\tau}dz \otimes d \bar{z}$);\\
(b)\ \ $g(z,w) =\log dist_\omega(z,w)+O(1), \,$ as $z\to w\,;$\\
(c)\ \ $g(w, w)=-\infty.$

Define further $c_{X_\tau}(z):=\exp\displaystyle \lim_{w \to z}(g(z,w)-\log dist_\omega(z,w))$ to be the modified logarithmic capacity, and then our research interest virtually focuses on the relations between $\pi K_{X_\tau}(z)$  and $c_{X_\tau}^2(z).$

Most strikingly, this modified green function\,(namely Arakelov-Green's function (see [10]),\,which we will explain later)\,has such a clear meaning in physics, pointed out by Ooguri in his recent famous lecture notes (see [11]). This may lead our results to some very practical applications.

\section{Proof of the main theorem}

Now we will make some detailed computations here to obtain the main theorem. From [12], we know the solution of the equation
$$\hspace{-1.5cm}   \frac{\partial ^2 g(z,w)}{\partial z \bar \partial z}=-\frac{\pi}{2}\cdot\frac{1}{ Im \tau}$$

\noindent is explicitly expressed as:
%$$\hspace{-1.0cm}   g(z,w)=\log \frac { \Vert \theta \Vert(z-w+\frac{1+\tau}{2};\tau)} {\Vert \eta \Vert(\tau)}.$$
\begin{equation}\label{1.1}
%\begin{split}
\hspace{-1.0cm}   g(z,w)=\log \frac { \Vert \theta \Vert(z-w+\frac{1+\tau}{2};\tau)} {\Vert \eta \Vert(\tau)}.
%\end{split}
 \end{equation}

\noindent where $ \Vert \theta \Vert(x+iy;\tau)=({Im \,\tau})^{\frac{1}{4}} \cdot \exp({-\pi y^2/{Im \,\tau}})\cdot |\theta(x+iy;\tau)|$ and $ \Vert \eta \Vert(\tau)=({Im \,\tau})^{\frac{1}{4}}\cdot |\eta (\tau)|$, and we know that $\theta(z;\tau):=\sum_{n=-\infty}^{\infty}\exp(\pi i n^2\tau+2\pi inz)\,\,$ is the theta function and $\eta (\tau):= q^{\frac{1}{12}} \cdot \prod_{n=1}^\infty (1-q^{2n})$ is the Dedekin-$\eta$ function, if we denote $\exp(\pi i \tau)$ by $q$.

Substituting these notifications into (1), we will get
$$g(z,w)=\log \frac { \exp({-\pi \big ({Im \, (z-w+\frac {\tau}{2})}\big)^2/{Im \,\tau}})\cdot |\theta(z-w+\frac{1+\tau}{2};\tau)|} {|q^{\frac{1}{12}} \cdot \prod_{n=1}^\infty (1-q^{2n})|}.$$

From [13, p.17], it holds that
$$\theta(z+\frac{1+\tau}{2};\tau)=\exp (-\frac{1}{4} \pi i \tau-\pi i (z+\frac{1}{2}))\sum_{n=-\infty}^{\infty}\exp(\pi i (n+\frac{1}{2})^2\tau+2\pi i(n+\frac{1}{2})(z+\frac{1}{2}))$$

Applying Jacobi triple product (see [14, Theorem 14.6]), since $q=\exp(\pi i \tau)$, we know $$g(z,w)=\log   \big ( e^{-\pi \big ({Im \, (z-w+\frac {\tau}{2})}\big)^2/{Im \,\tau}}  |2q^{\frac{1}{6}}\sin (\pi (z-w))\prod_{n=1}^\infty(1-2\cos(2\pi (z-w))q^{2n}+q^{4n}) \exp (-\frac{1}{4} \pi i \tau-\pi i (z-w+\frac{1}{2})) |\big).$$

By definition, it follows that

$$c_{X_\tau}(z)=\exp\displaystyle \lim_{w \to z}\Big(g(z,w)-\log dist_\omega(z,w)\Big)\hspace{5.3cm}$$
$$= ({ Im\, \tau})^{\frac{1}{2}}\cdot 2\pi\cdot\exp(-\frac{\pi}{6}Im\tau)\cdot |\,\prod_{n=1}^\infty(1-q^{2n})^2|.\hspace{2.6cm}$$

To acknowledge the relations between $\pi K_{X_\tau}$ and $c_{X_\tau}^2$, we define their quotient as a new function and
$$\hspace{-8.5cm}F(\tau):=\log \frac{\pi K_{X_\tau}}{c^2_{X_\tau}}$$
\begin{equation}
%\begin{split}
\hspace{-0.5cm}=-2\log (Im\, \tau)-\log( 4\pi)+\frac{\pi }{3} Im\tau-4\sum_{n=1}^\infty \log|1-q^{2n}|.
%\end{split}
\end{equation}

With the help of computers, we obtained that  $F(\tau)\geq -1.8251$ and therefore $\exp F(\tau)\geq  0.1612,\, \forall\, \tau \in \mathbb C\,(Im \tau>0).$ This means for any complex torus, the following inequality always holds:

$$\alpha\pi K_{X_\tau} \geq c^2_{X_\tau},\,\,\,\,\alpha \approx 6.2034.$$
And the $``="$ is attained when $Im \tau \approx 1.9192$ (Explicit programs for this numerical optimization are attached in the appendix).

%\hspace{2.5cm}\includegraphics[scale=0.45]{4.png}\\

\begin{figure}[htp]
  \centering
  \includegraphics[width=0.35\textwidth]{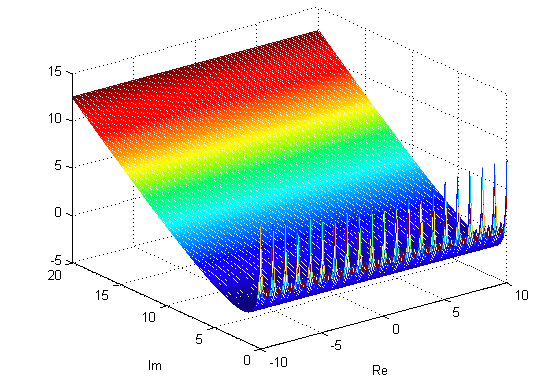}
   \caption{The 3D-graph of $F(\tau)$}
\end{figure}

Since the last term of $F(\tau)$ converges and tends to 0 as $Im \tau \to +\infty$, in particular, we have $$\lim_{Im \tau \to +\infty}\,\frac{\pi K_{X_\tau}}{c^2_{X_\tau}}=+\infty.$$

This result can be read as: For this special torus, the above $\alpha$ can be close to 0.

Now the main theorem and corollary are proved.

\section{Remarks on higher-genus cases}

In this section, we will give some remarks on the general cases for a compact Riemann surface with genus $g \geq 2$.

Let $X$ be a compact Riemann surface with genus $g \geq 2$ (more generally, hyperbolic in the sense of Kobayashi),
then there exists a conformal metric $\omega$ on $X$, obtained by descending the Poincare metric from its universal covering space to itself,\, s.t.  $Curv\,{\omega} \equiv -1.$ 

Consider the function $g(p,q):X\times X \to[-\infty,0),$ s.t. for each fixed $q \in X:$ \\
(a)\ \ $\Delta_{\omega} g(\cdot,q)=-1$  on $X\setminus \{q\}$, where $\Delta_{\omega}$ is the Laplacian with respect to the metric $\omega$;\\
(b)\ \ $\,g(p,q) = \log dist_\omega(p,q)+O(1), \,$ as $z\to w\,;$\\
(c)\ \ $g(w, w)=-\infty.$

The existence and uniqueness of such a function for general compact Riemann surfaces are proved by Arakelov (see [15]). And define further $c_X(p):=\exp\displaystyle \lim_{q \to p}\Big(g(p,q)- log(dist_\omega(p,q)) \Big)$ to be the logarithmic capacity with respect to $p$.

Let $H^{1,0}_{(2)}(X)$ be the $g-$dimensional Hilbert space of holomorphic 1-forms on $X$ such that $\big|\int_{X} f\wedge \bar f\,\big| <\infty\big\}.$ $H^{1,0}_{(2)}(X)$ naturally has an orthonormal basis $\big\{ \varphi_1,\cdots,\varphi_g\big\},$ and the Bergman kernel on the diagonal is of the (1,1)-form on $X,$ given by $$K_{X}(p)=\sum_{j=1}^g {\varphi_j(p)\wedge\overline {\varphi_j(p)}}.$$

So the remaining problem is to study the relations between $c^2_X(p)$ and $\pi K_X(p).$

\subsection*{Acknowledgements} 
\indent 
This work is supported by the National Natural Science Foundation of China (No.11031008, No.11171255) and the Nagoya University Gaku-Sei Project 2012. The author would like to express his sincere thanks to Professor Takeo Ohsawa for introducing him to this topic, giving many advices and reading the drafts. He would also like to thank M. Adachi and X. Liu for the helpful discussions. Last but not least, he would express his gratitude to J. Gustavo, T. Saikawa and N. Zheng for helping with the computer programming parts.

%\newpage

\subsection*{Appendix: \,programs on MATLAB}
\fontsize{11}{9}
\selectfont
To make this paper complete, we attach here the explicit programs, to get the numerical results concerning $F(\tau)$ defined in (2), running on MATLAB.\\

We first define a function called test(x,y,N) as below:\\

\hspace{0.5cm}\includegraphics[scale=0.5]{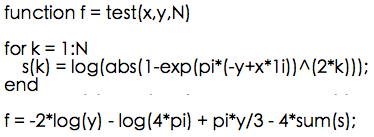}\\

Then we make the following computations:\\

\hspace{0.5cm}\includegraphics[scale=0.5]{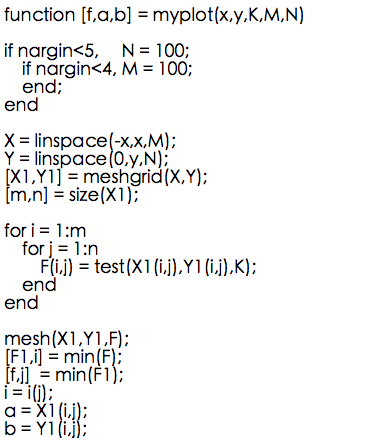}\\

Here $[-x, x]\times[0,y]$ forms the region where we plot the 3D-graph of $F(\tau)$. $f$ represents the minimal value on it for $F(\tau)$, achieved when $\tau=a+bi$. And $K,$ chosen to be large enough, denotes the total times of summation needed to get the desired precision. \\

To run this program, we may first choose the appropriate $x^*,y^*$ and $K^*$. Then it suffices to type into the command window:\\

$\,\,\,\,\,\,\,\,\,\,\,\,\,\,\,\gg\gg\,\,\,\,\,\,\,\,\,\,\,[f,a,b] = myplot(x^*,y^*,K^*)\,\,\,\,\,\,\,\,\,\,\,\ll\ll\hspace{5.8cm}$\\

\end{document}